\documentclass[12pt]{article}

\usepackage{amssymb,amsmath}  
\usepackage{bm}

\def\suml{\mathop{\sum}\limits}                                 
\def\maxl{\mathop{\max}\limits}                                 
\def\e{\varepsilon}                                             
\def\G{\Gamma}                                                  
\def\cd{\!\cdot}                                                
\DeclareMathSymbol{\ell}{\mathord}{letters}{96}                 
\def\id{\operatorname{id}}                                      
\def\od{\operatorname{od}}                                      
\def\FF{\mathop{\cal F}\nolimits}                               
\def\RR{\mathop{\cal R}\nolimits}                               
\def\rr{\mathop{\mathsf r}\nolimits}                            
\def\Fij{\FF^{\,i\to\bm{j}}}                                      
\def\Fto{\FF^{\,\to\bullet}}                                      
\def\cdc{,\ldots,}                                              
\def\1n{1\cdc n}                                                
\def\squareforqed{\hbox{\rlap{$\sqcap$}$\sqcup$}}               
\def\qed{\ifmmode\squareforqed\else{\unskip\nobreak\hfil\penalty50\hskip1em\null\nobreak\hfil\squareforqed 
\parfillskip=0pt\finalhyphendemerits=0\endgraf}\fi}             
\def\eq#1{\begin{equation}#1\end{equation}}                     
\def\eqs*#1{\begin{eqnarray*}#1\end{eqnarray*}}                 
\newtheorem{thm}{Theorem}{\bfseries}{\itshape}                  
\newtheorem{prop}{Proposition}{\bfseries}{\itshape}             
\newtheorem{corol}{Corollary}{\bfseries}{\itshape}              
\def\proof{{\noindent\bf Proof. }}                              
\providecommand{\url}[1]{#1}
\csname url@samestyle\endcsname

\title{A Graph Bottleneck Inequality}

\author{Pavel Chebotarev\\
       {\normalsize Institute of Control Sciences of the Russian Academy of Sciences}\\
       {\normalsize 65 Profsoyuznaya Street, Moscow 117997, Russia}\\
       {\small\tt chv@member.ams.org}
       }
\date{}

\begin{document}
\maketitle
\unitlength 1.50mm

\begin{abstract}
For a weighted directed multigraph, let $f_{ij}$ be the total weight of spanning converging forests that have vertex $i$ in a tree converging to~$j$. We prove that $f_{ij}\,f_{jk}=f_{ik}\,f_{jj}$ if and only if every directed path from $i$ to $k$ contains~$j$ (a~\emph{graph bottleneck equality}). Otherwise, $f_{ij}\,f_{jk}<f_{ik}\,f_{jj}$  (a~\emph{graph bottleneck inequality}). In a companion paper~\cite{Che08dist}
(P.~Chebotarev, A new family of graph distances, arXiv preprint math.CO/0810.2717, 2008. \url{http://arXiv.org/abs/0810.2717}. Submitted), this inequality underlies, by ensuring the triangle inequality, the construction of a new family of graph distances. This stems from the fact that the graph bottleneck inequality is a multiplicative counterpart of the triangle inequality for proximities.
\bigskip

\noindent{\em Keywords}:
 Spanning converging forest;
 Matrix forest theorem;
 Laplacian matrix
\bigskip

\noindent{\em AMS Classification}:
05C50, 
05C05, 
15A51
\end{abstract}

\section{Introduction}

Let $\G$ be a weighted directed multigraph with vertex set $V(\G)=\{\1n\}$, $n>1$. We assume that $\G$ has no loops.
For $i,j\in V(\G)$, let $n_{ij}\in\{0,1,\ldots\}$ be the number of arcs emanating from $i$ to $j$ in~$\G$; for every $p\in\{\1n_{ij}\}$, let $w_{ij}^p>0$ be the weight of the $p$th arc directed from $i$ to $j$ in~$\G$; let $w_{ij}=\sum_{p=1}^{n_{ij}}w_{ij}^p$ (if $n_{ij}=0$, we set $w_{ij}=0$) and $W=(w_{ij})_{n\times n}$. $W$ is the \emph{matrix of total arc weights}. The \emph{outdegree\/} and \emph{indegree\/} of vertex $i$ are $\od(i)=\sum_{j=1}^nn_{ij}$ and $\id(i)=\sum_{j=1}^nn_{ji}$, respectively.
\smallskip

A \emph{converging tree\/} is a weakly connected 
weighted digraph in which one vertex, called the {\it root}, has outdegree zero and the remaining vertices have outdegree one.
A~{\it converging forest\/} is a weighted digraph all of whose weakly connected components are converging trees. The roots of these trees are referred to as the roots of the converging forest. A spanning converging forest of $\G$ is called an {\em in-forest\/} of~$\G$.
\smallskip

By the weight of a weighted digraph $H$, $w(H)$, we mean the product of the weights of all its arcs. If $H$ has no arcs, then $w(H)=1$. The weight of a set $S$ of digraphs, $w(S)$, is the sum of the weights of the digraphs belonging to~$S$; the weight of the empty set is zero. If the weights of all arcs are unity, i.\,e., the graphs in $S$ are actually unweighted, then $w(S)$ reduces to the cardinality of~$S$.
\smallskip

For a fixed $\G$, by $\Fto$ and $\Fij$ we denote the set of all in-forests of $\G$ and the set of all in-forests of $\G$ that have vertex $i$ belonging to a tree rooted at~$j$, respectively. Let {$f=w(\Fto)$} and
\eq{
\label{e_fij}
f_{ij}=w(\Fij),\quad i,j\in V(\G);
}
by $F$ we denote the matrix with entries $f_{ij}$: $F=(f_{ij})_{n\times n}$. $F$ is called the \emph{matrix of in-forests of~\,$\G$}.
\smallskip

Let $L=(\ell_{ij})$ be the Laplacian matrix of $\G$, i.\,e.,
$$
\ell_{ij}=
        \begin{cases}
        -w_{ij},               &j\ne i,\\
         \suml_{k\ne i}w_{ik}, &j  = i.
        \end{cases}
$$

Consider the matrix 
\eq{
\label{e_Q}
Q=(q_{ij})=(I+L)^{-1}.
}
By the matrix forest theorem\footnote{Versions of this theorem for undirected (multi)graphs can be found in~\cite{CheSha95,Merris97}.} \cite{CheSha97,CheAga02ap}, for any weighted digraph $\G$, $Q$ does exist and
\eq{
\label{e_mft}
q_{ij}=\frac{f_{ij}}{f},\quad i,j=\1n.
}
Therefore $F=fQ=f\cd(I+L)^{-1}$. The matrix $Q$ can be considered as a proximity (similarity) matrix of~$\G$ \cite{CheSha97,Che08DAM}.
\smallskip

In Section~\ref{s_main}, we present the \emph{graph bottleneck inequality\/} involving the $f_{ij}$'s and a necessary and sufficient condition of its reduction to equality.

\section{A graph bottleneck inequality and a graph\\ bottleneck equality}
\label{s_main}

\begin{thm}
\label{th_main}
Let $\G$ be a weighted directed multigraph and let the values $f_{ij}$ be defined by~\eqref{e_fij}. Then for every $i,j,k\in V(\G),$
\eq{
\label{e_ineq}
f_{ij}\,f_{jk}\le f_{ik}\,f_{jj}.
}
Moreover$,$
\eq{
\label{e_eq}
f_{ij}\,f_{jk}=f_{ik}\,f_{jj}
}
if and only if every directed path from $i$ to $k$ contains~$j$.
\end{thm}

Since \eqref{e_ineq} reduces to \eqref{e_eq} when $j$ is a kind of bottleneck in $\G$, \eqref{e_eq} is called a \emph{graph bottleneck equality}; by the same reason, \eqref{e_ineq} is referred to as a \emph{graph bottleneck inequality}. It is readily seen that the graph bottleneck inequality is a multiplicative counterpart of the triangle inequality for proximities (see, e.\,g.,~\cite{CheSha97}).
\smallskip

It turns out that it is not easy to construct a direct bijective proof to Theorem~\ref{th_main}.
We present a different proof; 
it requires some additional notation and two propositions given below.
\smallskip

For a fixed multidigraph $\G$,
let us choose an arbitrary $\e>0$ such that $0\le\e\cd\!\maxl_{1\le i\le n}\ell_{ii}<1$.
It is easy to verify that the matrix
\eq{
\label{e_Pm}
P=(p_{ij})=I-\e L
}
is row stochastic: $0\le p_{ij}\le1$ and $\sum_{k=1}^n p_{ik}=1$, $\;i,j=\1n$.
\smallskip\smallskip

Denote by $\G^{\circlearrowright}$ the weighted multidigraph with loops whose matrix $W(\G^{\circlearrowright})$ of total arc weights is $(1+\e)^{-1}P$. More specifically, every vertex $i$ of $\G^{\circlearrowright}$ has a loop with weight $(1+\e)^{-1}p_{ii}$; the remaining arcs of $\G^{\circlearrowright}$ are the same as in $\G$, their weights being the corresponding weights in~$\G$ multiplied by $(1+\e)^{-1}\e$.

Recall that a $v_0\to v_k$ \emph{route\/} in a multidigraph with loops is an alternating sequences of vertices and arcs $v_0,x_1,v_1\cdc x_k,v_k$ where each arc $x_i$ is $(v_{i-1},v_i)$.
The \emph{length\/} of a route is the number $k$ of its arcs (including loops). The \emph{weight\/} of a route is the product of the weights of all its arcs.
We assume that for every vertex $i$, there is a unique route of length $0$ from $i$ to~$i$, the weight of this route being~1. The \emph{weight of a set of routes\/} is the total weight of the routes the set contains. 
\smallskip

Let $r_{ij}$ be the weight of the set $\RR^{ij}$ of all $i\to j$ routes in $\G^{\circlearrowright}$, provided that this weight is finite (this reservation is essential because the set of $i\to j$ routes is infinite whenever $j$ is reachable from~$i$). $R=(r_{ij})_{n\times n}$ will denote the \emph{matrix of the total weights of routes}. 

\begin{prop}
\label{p_RoutFo}
For every weighted multidigraph $\G$ and every $\e>0$ such that $0\le\e\cd\!\maxl_{1\le i\le n}\ell_{ii}<1,$ the matrix $R$ of the total weights of routes in $\G^{\circlearrowright}$ exists and it is proportional to the matrix $F$\! of in-forests of~$\G$.
\end{prop}

\proof
Observe that for every $k=0,1,2,\ldots,$ the matrix of total weights of $k$-length routes in $\G^{\circlearrowright}$ is $((1+\e)^{-1}P)^k$. Therefore the matrix $R$, whenever it exists, can be expressed as follows:
\eq{
\label{e_RfrP}
R=\suml_{k=0}^\infty((1+\e)^{-1}P)^k.
}
Since the spectral radius of $P$ is 1 and $0<(1+\e)^{-1}<1$, the sum in \eqref{e_RfrP} does exist\footnote{On counting routes see~\cite{Kasteleyn67}.  Related finite topological representations that involved paths were obtained in~\cite{Ponstein66}. For a connection with matroid theory see, e.\,g.,~\cite{Schrijver78}.}, 
therefore \eqref{e_RfrP}, \eqref{e_Pm}, \eqref{e_Q}, and \eqref{e_mft} imply
\eqs*{
\label{e_RfrP1}
R
&=&(I-(1+\e)^{-1}P)^{-1}
 = \left(I-(1+\e)^{-1}\!\left(I-\e L\right)\!\right)^{-1}\\
&=&\left(\frac{\e}{1+\e}(I+L)\!\right)^{-1}
 = \left(1+\e^{-1}\right)Q
 = \left(1+\e^{-1}\right)f^{-1}F,
}
which completes the proof.\qed

\begin{prop}
\label{p_Routes}
For any weighted multidigraph with loops and any vertices $i,j,$ and $k,$ if the total weights of routes $r_{ij},$ $r\!_{jj},$ $r\!_{jk},$ and $r_{ik}$ are finite$,$ then
\eq{
\label{e_ineq1}
r_{ij}\,r\!_{jk}\le r_{ik}\,r\!_{jj}.
}
Moreover$,$
\eq{
\label{e_eq1}
r_{ij}\,r\!_{jk}=r_{ik}\,r\!_{jj}
}
if and only if every directed path from $i$ to $k$ contains~$j$.
\end{prop}

\proof
Suppose that the total weights of routes $r_{ij},$ $r\!_{jj},$ $r\!_{jk},$ and $r_{ik}$ are finite.
Let $\RR^{ij(1)}$ be the set of all $i\to j$ routes that contain only one appearance of~$j$. Let $r_{ij(1)}=w(\RR^{ij(1)})$. Then every $i\to j$ route $\rr^{ij}\in \RR^{ij}$ can be uniquely decomposed into a route $\rr^{ij(1)}\in\RR^{ij(1)}$ and a route (possibly, of length~0) $\rr^{jj}\in\RR^{jj}$. And vice versa, linking an arbitrary route $\rr^{ij(1)}\in\RR^{ij(1)}$ with an arbitrary $\rr^{jj}\in\RR^{jj}$ results in a certain route $\rr^{ij}\in\RR^{ij}$. This determines a natural bijection between $\RR^{ij}$ and $\RR^{ij(1)}\times\RR^{jj}$. Therefore
\eq{
\label{e_rij}
r_{ij}=r_{ij(1)}\,r\!_{jj}.
}

Let $\RR^{ijk}$ and $\RR^{i\bar\jmath k}$ be the sets of all $i\to k$ routes that contain and do not contain $j$, respectively. Then $\RR^{ik}=\RR^{ijk}\cup\RR^{i\bar\jmath k}$ and $\RR^{ijk}\cap\RR^{i\bar\jmath k}=\varnothing$, consequently,
\eq{
\label{e_rik}
r_{ik}=r_{ijk}+r_{i\bar\jmath k},
}
where $r_{ijk}=w(\RR^{ijk})$ and $r_{i\bar\jmath k}=w(\RR^{i\bar\jmath k})$.
\smallskip

Furthermore, by the argument similar to that justifying \eqref{e_rij} one obtains
\eq{
\label{e_rijk}
r_{ijk}=r_{ij(1)}\,r\!_{jk}.
}

Combining \eqref{e_rik}, \eqref{e_rijk}, and \eqref{e_rij} yields
\eqs*{
\label{e_rjjik}
 r_{ik}\,r\!_{jj}
&=&(r_{ijk}+r_{i\bar\jmath k})\,r\!_{jj}
 = r_{ij(1)}\,r\!_{jk}\,r\!_{jj}+r_{i\bar\jmath k}\,r\!_{jj}\\
&=&r_{ij}\,r\!_{jk}+r_{i\bar\jmath k}\,r\!_{jj}\geq r_{ij}\,r\!_{jk},
}
with the equality if and only if $\RR^{i\bar\jmath k}=\varnothing$.
\qed
\bigskip

{\noindent\bf Proof of Theorem~\ref{th_main}.}
Theorem~\ref{th_main} follows immediately by combining Propositions~\ref{p_RoutFo} and~\ref{p_Routes}.
\qed
\bigskip

Finally, consider the graph bottleneck inequality and the graph bottleneck equality for undirected graphs.

\begin{corol}[{\bf to Theorem~\ref{th_main}}]
\label{co_main}
Let $G$ be a weighted undirected multigraph and let $f_{ij},$ $i,j\in V(G),$ be the total weight of all spanning rooted forests of $G$ that have vertex $i$ belonging to a tree rooted at~$j$. Then for every $i,j,k\in V(G),$
\eq{
\label{e_uineq}
f_{ij}\,f_{jk}\le f_{ik}\,f_{jj}.
}
Moreover$,$
\eq{
\label{e_ueq}
f_{ij}\,f_{jk}=f_{ik}\,f_{jj}
}
if and only if every path from $i$ to $k$ contains~$j$.
\end{corol}

\proof
Consider the weighted multidigraph $\G$ obtained from $G$ by replacing every edge by two opposite arcs carrying the weight of that edge. Then, by the matrix forest theorems for weighted and unweighted graphs, $f_{ij}(G)=f_{ij}(\G),$ $i,j\in V(G)$. Observe that for every $i,j,k\in V(G),$ every path from $i$ to $k$ contains~$j\/$ if and only if every directed path in $\G$ from $i$ to $k$ contains~$j$. Therefore, by virtue of Theorem~\ref{th_main}, inequality \eqref{e_uineq} follows for $G$; moreover, equality \eqref{e_ueq} holds true if and only if every path in $G$ from $i$ to $k$ contains~$j$.
\qed
\medskip

In a companion paper \cite{Che08dist}, the graph bottleneck inequality for undirected graphs is used to ensure the triangle inequality for a new parametric family $\{d_\alpha(\cdot,\cdot)\}$ of graph distances. In turn, the bottleneck equality provides a necessary and sufficient condition under which the triangle inequality $d_\alpha(i,j)+d_\alpha(j,k)\ge d_\alpha(i,k)$ for a triple $i,j,k$ of graph vertices reduces to equality.

\end{document}